\documentclass[
  letterpaper,
  10 pt,
  conference
]{ieeeconf}
\IEEEoverridecommandlockouts
\makeatletter
\let\NAT@parse\undefined
\makeatother

\usepackage{amsmath,amssymb}
\usepackage{xcolor}
\usepackage{supertabular,booktabs,array}
\usepackage{flushend}
\usepackage{graphicx}
\usepackage{float}
\usepackage[mathscr]{eucal}
\usepackage{hyperref,cleveref}

\title{Primal-Dual iLQR}
\author{
João Sousa-Pinto (Apple) \\ \href{mailto:joaospinto@gmail.com}{joaospinto@gmail.com}
\and Dominique Orban (Polytechnique Montréal) \\ \href{mailto:dominique.orban@polymtl.ca}{dominique.orban@polymtl.ca}
}
\date{December 2023}

\begin{document}
\maketitle

\section{Abstract}\label{abstract}

We introduce a new algorithm for solving unconstrained discrete-time
optimal control problems. Our method follows a direct multiple shooting
approach, and consists of applying the SQP method together with an
\(\ell_2\) augmented Lagrangian primal-dual merit function.
We use the LQR algorithm to efficiently solve the primal-dual Newton-KKT system.
As our algorithm is a specialization of NPSQP~\cite{ref-NPSQP}, it
inherits its generic properties, including global convergence, fast
local convergence, and the lack of need for second order corrections or dimension expansions,
improving on existing direct multiple shooting approaches such as
acados~\cite{ref-ACADOS}, ALTRO~\cite{ref-ALTRO}, GNMS~\cite{ref-GNMS},
FATROP~\cite{ref-FATROP}, and FDDP~\cite{ref-CROCODDYL}.
The solutions of the LQR-shaped subproblems posed by our algorithm can be
be parallelized to run in time logarithmic in the number of stages, states, and controls.
Moreover, as our method avoids sequential rollouts of the nonlinear dynamics,
it can run in $O(1)$ parallel time per line search iteration.
Therefore, this paper provides a practical, theoretically sound, and highly parallelizable
(for example, with a GPU) method for solving nonlinear discrete-time optimal control problems.
An open-source JAX implementation of this algorithm can be found on
\href{https://github.com/joaospinto/primal_dual_ilqr}
{GitHub ($joaospinto/primal \textunderscore dual \textunderscore ilqr$)}.

\section{Introduction}\label{introduction}

\subsection{Unconstrained Discrete-Time Optimal Control
Problems}\label{unconstrained-discrete-time-optimal-control-problems}

Unconstrained discrete-time optimal control problems are optimization
problems of the form

\begin{math}
\begin{aligned}
&\min_{x_0, u_0 \ldots, x_N} &\sum \limits_{i=0}^{N-1} g_i(x_i, u_i) + g_N(x_N) \\
&\mbox{s.t.} & x_0 = s_0 \\
& & x_{i+1} = f_i(x_i, u_i), \forall i \in \lbrace 0, \ldots, N-1 \rbrace.
\end{aligned}
\end{math}

Such optimization problems are ubiquitous in the fields of motion
planning and controls.

\subsection{Related Work}\label{related-work}

Naïf applications of generic optimization methods to unconstrained
discrete-time optimal control problems would require
\(O(N^3 (n + m)^3)\) operations for solving the linear systems posed at
each iteration, where \(N, n, m\) are respectively the number of stages,
states, and controls of the problem. The first algorithm to improve on
this was DDP~\cite{ref-DDP}, achieving complexity \(O(N (n + m)^3)\). Its
convergence properties were studied in~\cite{ref-ConvergenceDDP}.
iLQR~\cite{ref-ILQR} can be seen as a version of DDP with coarser
second derivative information (in particular, disregarding the Hessians
of the dynamics), trading off local quadratic convergence for some added
ease of implementation. Stagewise Newton~\cite{ref-StagewiseNewton}
exploits the fact that the computation of a Newton step on the
``eliminated problem'' (obtained by eliminating the variables \(x_i\) in
decreasing order of \(i\) by plugging in the dynamics into the costs)
can be reduced to solving an LQR problem, thereby achieving the same
computational complexity as DDP.

DDP, iLQR, and Stagewise Newton are all single shooting methods. This
makes them hard to warm start, as the state trajectory cannot be chosen
independently of the control trajectory. This can result in more
iterations being required for these methods to converge, as well as
higher likelihoods of convergence to undesired local optima. Moreover,
they require that the dynamics be evaluated sequentially, making them
harder to parallelize.

\cite{ref-GNMS} introduces a multiple shooting equivalent of
iLQR, called GNMS, as well as an algorithm that combines single and
multiple shooting, called iLQR-GNMS, with similar computational
complexity to the methods discussed above. There are no available
convergence results for iLQR-GNMS; at a minimum, a line search
procedure or a filter method would have to be added for these algorithms
to achieve global convergence. \cite{ref-CROCODDYL} introduces an
alternative multiple shooting equivalent of iLQR, called FDDP, which
adds a curvilinear line search that smoothly transitions between using
the nonlinear dynamics in the forward pass for larger step sizes and
using the affine dynamics for smaller step sizes. There are also no
available convergence results for FDDP. In both cases, as second order
derivatives of the dynamics are not considered, there cannot be local
quadratic convergence.

Much work has also been done in handling discrete-time optimal control
problems with stagewise constraints. ALTRO~\cite{ref-ALTRO},
FATROP~\cite{ref-FATROP}, and acados~\cite{ref-ACADOS}
are likely the most competitive open source software
packages for this application at the time of writing. ALTRO handles
constraints via the augmented Lagrangian method, while FATROP uses an
interior point method and acados uses an SQP approach.
ForcesNLP~\cite{ref-FORCESNLP} is likely the most competitive competing
commercial software package at the time of writing, and also relies on an
interior point method.

At a high level and modulo technicalities, FATROP can be seen as a
specialization of IPOPT~\cite{ref-IPOPT} with a linear system
solver that handles LQR problems with extra stagewise affine constraints
in time \(O(N (n + m)^3)\). Inequality constraints are handled by
eliminating some of the block-rows of the modified KKT system, as
originally done in~\cite{ref-InteriorPointMPC}. The modified KKT
system will only have the extra stagewise affine constraints when the
problem being solved has other stagewise equality constraints besides
the dynamics. A similar algorithm that does not handle other general
equality constraints had previously been presented in~\cite{ref-InteriorPointDDP}.
As compared to the specialization of FATROP to
problems without stagewise equality constraints, our method has the
advantage of not requiring second order corrections for fast local
convergence. ForcesNLP also follows an interior point approach, but
relies on a Cholesky factorizations instead of LQR decompositions for
solving its linear systems.

On the other hand, the augmented Lagrangian technique employed by ALTRO
yields linear systems that can be solved via the LQR algorithm, even in
the presence of stagewise equality constraints. However, ALTRO has a
slightly unorthodox way of doing multiple shooting: it modifies the
dynamics \(x_{n+1} = f_n(x_n, u_n)\) into
\(x_{n+1} = f_n(x_n, u_n) + e_n\) and adds the constraints \(e_n = 0\).
This makes it possible to independently warm start the state trajectory
(while breaking the non-dynamic \(e_n = 0\) constraint) while still
using a single shooting algorithm. This approach has the disadvantage of
substantially expanding the control dimension.

Finally, the SQP approach followed by acados poses
inequality-constrained quadratic subproblems, which are then solved by
HPIPM~\cite{ref-HPIPM}, an efficient interior point QP solver
that exploits the discrete-time optimal control problem structure. While
it does not currently support stagewise equality constraints other than
the dynamics, it could be extended to do so using a similar approach
to~\cite{ref-FATROP}. Moreover, due to not employing a globalization
strategy, global convergence is not guaranteed. When the functions
defining the problem are expensive to evaluate, SQP methods become
time-competitive, as the higher number of linear system solves will be
compensated by the lower number of evaluations of the nonlinear problem.

\subsection{Contributions}\label{contributions}

This paper introduces the Primal-Dual iLQR algorithm, which consists of a
specialization of the generic NPSQP algorithm to the case of unconstrained
discrete-time optimal control problems.

Our algorithm explores the sparsity structure of the Newton-KKT systems of
unconstrained discrete-time optimal control problems, resulting in only
\(O(N (n + m)^3)\) operations for solving the linear systems posed at
each iteration, where \(N, n, m\) are respectively the number of stages,
states, and controls of the problem.

Moreover, we show that the primal-dual Newton-KKT systems we pose can be solved
in $O(\log(N) \log(n) + \log(m))$ parallel time complexity, using a minor variation of~\cite{ref-AS-LQR},
and that the remaining parts of our method have constant parallel time complexity.

Our algorithm improves on earlier direct multiple shooting methods,
specifically by guaranteeing global convergence and
not impeding local superlinear convergence (even without second order corrections),
under the assumptions described in~\cite{ref-NPSQP}.

\section{General Optimization
Background}\label{general-optimization-background}

\subsection{First-order Optimality Conditions and KKT
Systems}\label{first-order-optimality-conditions-and-kkt-systems}

Given a constrained optimization problem of the form \begin{align}
&\min_{x} &g(x) \\
&\mbox{s.t.} &c(x)=0,
\end{align} for any local optimum \(x\) meeting the standard linear
independence constraint qualification (LICQ) constraints (i.e., for
which the rows of \(J(c)(x)\) are linearly independent), there exists a
vector \(\lambda\) (called the Lagrange/KKT multiplier) for which
\begin{align}
\begin{cases}
\nabla g(x) + \nabla c(x) \lambda &= 0 \\
c(x) &= 0.
\end{cases}
\end{align} In other words, \((x, \lambda)\) is a critical point (but
not necessarily a minimizer) of the Lagrangian function
\(\mathcal{L}(x, \lambda)\) defined as \(g(x) + \lambda^T c(x)\).

\subsection{Sequential Quadratic
Programming}\label{sequential-quadratic-programming}

The sequential quadratic programming (SQP) method for
equality-constrained non-linear optimization problems consists of
applying Newton's method for finding zeros of \(\nabla \mathcal{L}\),
typically in combination with a line search mechanism, in order to
ensure global convergence.

Below, let \(A = J(c)(x)\), \(l = \nabla_x \mathcal{L}(x, \lambda)\),
and \(d = c(x)\). Moreover, let \(Q\) be a positive definite
approximation of \(\nabla_{xx}^2 \mathcal{L}(x, \lambda)\). Such an
approximation can be constructed, for example, via regularization.
Alternatively, the terms involving \(\nabla_{xx}^2 c(x)\) can be dropped
entirely, and \(Q\) would consist purely of a positive definite
approximation of \(\nabla_{xx} g(x)\). The latter approach, of course,
would result in local quadratic convergence being lost. For simplicity
of notation, we henceforth drop the dependencies of these terms in
\(x\).

Each Newton step for finding a zero of \(h(z)\) is obtained by solving
\(h'(z) \Delta z = -h(z)\). If \(h = \nabla \mathcal{L}\), this
corresponds to solving \begin{align}
\begin{pmatrix}
\label{ref-newton-kkt}
Q & A^T \\
A & 0
\end{pmatrix}
\begin{pmatrix} \Delta x \\ \Delta \lambda \end{pmatrix} =
-\begin{pmatrix} l \\ d \end{pmatrix}.
\end{align}

It can be shown (by a simple application of the first-order optimality
conditions mentioned above) that \(\Delta x\) is the solution of
\begin{align}
&\min_p &\frac{1}{2} p^T Q p + l^T p \\
&\mbox{s.t.} &Ap + d = 0.
\end{align} Similarly, it can also be shown that \(\Delta \lambda\) is
the corresponding Lagrange multiplier. As long as \(Q\) is positive
definite, this cost function is bounded from below. Depending on the
sparsity pattern of this system, weaker requirements on \(Q\) might also
guarantee this.

Convergence is established when \(\Delta x\) and \(\Delta \lambda\) are
sufficiently close to \(0\). The solution is feasible as long as
\(d = 0\), but convergence to locally infeasible solutions is possible.

\subsection{Merit Functions}\label{merit-functions}

Merit functions offer a measure of progress towards the final solution
and are typically used in combination with a line search procedure. Note
that it would be inadequate to simply apply line search on the original
objective, as SQP may (and often does) provide a search direction that
locally increases that objective (usually trading that off for a
reduction in the constraint violations). Filter methods, used in
~\cite{ref-IPOPT}, provide an alternative mechanism for
measuring progress and determining whether the candidate step should be
accepted, but this paper will not be using that approach.

Given a penalty parameter \(\rho > 0\), we define the \(\ell_2\)
merit function \begin{align}
m_\rho(x, \lambda) = \mathcal{L}(x, \lambda) + \frac{\rho}{2} \lVert c(x) \rVert^2.
\end{align}

Below, \(\Delta x\) and \(\Delta \lambda\) correspond to SQP step
defined above; moreover, $D(\cdot; \cdot)$ is used to represent the
directional derivative operator. Noting that
\(\nabla_\lambda \mathcal{L}(x, \lambda) = c(x) = d\), it holds that
\begin{equation}
\begin{aligned}
&D(m_\rho; \begin{pmatrix} \Delta x^T & 0 \end{pmatrix}) \\
&= \Delta x^T \nabla_x m_\rho(x, \lambda) \\
&= \Delta x^T \nabla_x \mathcal{L}(x, \lambda) + \rho \Delta x^T J(c)(x)^T c(x) \\
&= -\Delta x^T Q \Delta x - \Delta x^T A^T \Delta \lambda + \rho \Delta x^T A^T d \\
&= -\Delta x^T Q \Delta x + \nabla_\lambda \mathcal{L}(x, \lambda)^T \Delta \lambda
- \rho \nabla_\lambda \mathcal{L}(x, \lambda)^T d \\
&= -\Delta x^T Q \Delta x + d^T \Delta \lambda - \rho \lVert d \rVert^2
\end{aligned}
\end{equation}
and that
\begin{equation}
\begin{aligned}
D(m_\rho; \begin{pmatrix} 0 & \Delta \lambda^T \end{pmatrix})
&= \nabla_\lambda \mathcal{L}(x, \lambda)^T \Delta \lambda = d^T \Delta \lambda.
\end{aligned}
\end{equation}
Therefore:
\begin{equation}
\begin{aligned}
&D(m_\rho; \begin{pmatrix} \Delta x^T & \Delta \lambda^T \end{pmatrix}) \\
= &-\Delta x^T Q \Delta x + 2 d^T \Delta \lambda - \rho \lVert d \rVert^2.
\end{aligned}
\end{equation}

Note that, as \(Q\) is positive definite, \(-\Delta x^T Q \Delta x < 0\)
unless \(\Delta x = 0\). When \(d\) is not acceptably close to \(0\), we
suggest taking \begin{align}
\rho = \frac{2 \lVert \Delta \lambda \rVert}{\lVert d \rVert},
\end{align} which ensures that
\(D(m_\rho; \begin{pmatrix} \Delta x^T & \Delta \lambda^T \end{pmatrix}) < 0\)
(due to the Cauchy-Schwarz inequality). Otherwise, we suggest taking
\(\rho = 0.01\).

Note that~\cite{ref-NPSQP} suggests instead never decreasing \(\rho\),
and whenever an increase is required, increasing it to twice the minimal
value required for making \begin{align}
D(m_\rho; \begin{pmatrix} \Delta x^T & \Delta \lambda^T \end{pmatrix}) < 0.
\end{align}

\subsection{Line Search}\label{line-search}

Having computed \(\Delta x\), \(\Delta \lambda\), and \(\rho\), we wish
to compute a step size \(\alpha\) that results in an acceptable decrease
in our merit function \(m_\rho\). Typically, this is achieved by
performing a backtracking line search. \(\alpha\) starts as \(1\) and
decreases by a constant multiplicative factor (often \(0.5\)) every
iteration, until
\begin{equation}
\begin{aligned}
&m_\rho(x + \alpha \Delta x, \lambda + \alpha \Delta \lambda) < \\
&m_\rho(x, \lambda) + k \alpha D(m_\rho; \begin{pmatrix} \Delta x^T & \Delta \lambda^T \end{pmatrix})
\end{aligned}
\end{equation}
is satisfied (where \(k \in (0, 1)\) is called the Armijo
factor; typically, \(10^{-4}\) is used). This is called the Armijo
condition. This process is guaranteed to terminate as long as \(g(x)\)
and \(c(x)\) are differentiable at \(x\) and
\(D(m_\rho; \begin{pmatrix} \Delta x^T & \Delta \lambda^T \end{pmatrix}) < 0\).

When \(\Delta x = 0\) and \(d = 0\) but \(\Delta \lambda \neq 0\), we
take a full step without conducting a line search.

Note that NPSQP~\cite{ref-NPSQP} uses a slightly different line search
method. Specifically, it enforces both Wolfe conditions, not only the
Armijo condition.

\section{Algorithm Derivation}\label{algorithm-derivation}

In this section, we specialize the methods described above to the case
of unconstrained discrete-time optimal control problems as needed.

We start by defining \begin{align}
d =
\begin{pmatrix}
d_0 \\
\cdots \\
d_N
\end{pmatrix} =
c(x) =
\begin{pmatrix}
s_0 - x_0 \\
f(x_0, u_0) - x_1 \\
\cdots \\
f(x_{N - 1}, u_{N - 1}) - x_N
\end{pmatrix}.
\end{align} The corresponding Lagrange/KKT multipliers are
\(\lambda_0, \ldots, \lambda_N\) respectively. Letting
\(x = (x_0, u_0, \ldots, x_{N - 1}, u_{N - 1}, x_N)\),
\(A_i = J_x(f_i)(x_i, u_i)\) and \(B_i = J_u(f_i)(x_i, u_i)\),
\begin{equation}
\begin{aligned}
& A = J(c)(x) = \\
& \begin{pmatrix}
-I   &  0    &  0    &  0    &  0       &  \cdots  &  0          &  0          &  0  \\
A_0  &  B_0  &  -I   &  0    &  0       &  \cdots  &  0          &  0          &  0  \\
0    &  0    &  A_1  &  B_1  &  -I      &  \cdots  &  0          &  0          &  0  \\
0    &  0    &  0    &  0    &  \ddots  &  \ddots  &  \ddots     &  0          &  0  \\
0    &  0    &  0    &  0    &  0       &  \cdots  &  A_{N - 1}  &  B_{N - 1}  &  -I
\end{pmatrix}.
\end{aligned}
\end{equation}

Note that the LICQ conditions are always met for matrices of this form,
due to the presence of the \(-I\) blocks, independently of the \(A_i\)
and \(B_i\).

The SQP constraint \(Ap + d = 0\) becomes \begin{align}
\begin{cases}
\Delta x_0 = d_0 \\
\Delta x_{i + 1} = A_i \Delta x_i + B_i \Delta u_i + d_i,
\forall i \in \lbrace 0, \ldots, N - 1 \rbrace.
\end{cases}
\end{align} Note that if \(x_0\) is warm started as \(s_0\) then
\(\Delta x_0\) will always be \(0\).

Similarly, noting that
\begin{equation}
\begin{aligned}
\mathcal{L}(x, \lambda) =
&\lambda_0^T (s_0 - x_0) + \sum\limits_{i=0}^{N - 1} g_i(x_i, u_i) + \\
& \lambda_{i + 1}^T (f_i(x_i, u_i) - x_{i + 1})
+ g_N(x_N),
\end{aligned}
\end{equation}
it holds that \begin{align}
l = \nabla_x \mathcal{L}(x, \lambda) =
\begin{pmatrix}
l_0 \\
\vdots \\
l_N
\end{pmatrix}
\end{align} where
\begin{equation}
\begin{aligned}
\begin{cases}
l_i = \nabla g_i(x_i, u_i) + \nabla f_i(x_i, u_i) \lambda_{i + 1} - \lambda_i, \\
\forall i \in \lbrace 0, \ldots, N - 1 \rbrace \\
l_N = \nabla g_N(x_N) - \lambda_N
\end{cases}
\end{aligned}
\end{equation}
and that \begin{align}
Q = \nabla_{xx} \mathcal{L}(x, \lambda) =
\begin{pmatrix}
Q_0 & & \\
& \ddots \\
& & Q_N
\end{pmatrix},
\end{align} where \begin{align}
\begin{cases}
Q_i = \nabla^2 g_i(x_i, u_i) + \nabla^2 f_i(x_i, u_i) \lambda_{i + 1},
\forall i \in \lbrace 0, \ldots, N - 1 \rbrace \\
Q_N = \nabla^2 g_N(x_N).
\end{cases}
\end{align} The SQP cost \(\frac{1}{2} p^T Q p + l^T p\) is simply
\(\sum \limits_{i=0}^N \frac{1}{2} p_i^T Q_i p_i + l_i^T p_i\), where
\(p_i = (\Delta x_i, \Delta u_i)\) for
\(i \in \lbrace 0, \ldots, N - 1 \rbrace\) and \(p_N = \Delta x_N\).

If exact Hessians are used, the matrices \(Q_i\) need not be positive
semi-definite. However, we require that positive definite approximations
be used instead. It would also be acceptable to only require positive
semi-definiteness, as long as the \(\nabla_{u_i u_i}^2\) components of
the \(Q_i\) be positive definite. However, this would require checking
that \(\Delta x^T Q \Delta x > 0\) at the end of the LQR solve. If this
condition is not met, positive definite approximations would have to be
used. Requiring positive definiteness in the first place avoids this
complication.

In this case, the the SQP problem is a primal-dual LQR problem,
which can be efficiently solved.

In order to ensure that the matrices \(Q_i\) are positive semi-definite,
it may be helpful to maintain regularization parameters \(\mu_i\) and
always use \(Q_i + \mu_i I\) instead of \(Q_i\). If, during the LQR
process, we detect that \(Q_i + \mu_i I\) is not positive definite,
\(\mu_i\) can be updated by multiplying it by an updated factor
\(r > 1\). When increasing \(\mu_i\) is not required, we can instead
update \(\mu_i\) by dividing it by \(r\). Minimum and maximum
regularization parameters may also be established, to allow, when
possible, \(\mu_i\) to eventually be set to \(0\) if \(Q_i\) is
consistently positive definite, as well as to prevent exploring
unreasonably high regularization parameters.
Another option would be to regularize the $Q_i$ by performing an
explicit eigenvalue decomposition and removing negative (or non-positive,
when positive definiteness is required) eigenvalues.

\section{Sequential Primal LQR Overview}\label{sequential-lqr}

In this section, we will go over the conventional sequential algorithm for solving primal LQR problems.
In the interest of not deviating from standard notation, variable names
may henceforth not match earlier parts of this paper.
A linear-quadratic regulator (LQR) problem is an unconstrained discrete-time optimal control problem
where the costs are quadratic functions and the dynamics are affine functions.
Specifically, they are optimization problems of the form
\begin{equation}
\begin{aligned}
&\min_{x_0, u_0 \ldots, x_N} &\sum \limits_{i=0}^{N-1}
  \left(
  \frac{1}{2} x_{i}^{T} Q_{i} x_{i} + q_{i}^{T} x_{i} +
  \frac{1}{2} u_{i}^{T} R_{i} u_{i} + r_{i}^{T} u_{i} \right. \\
  && \left. \vphantom{\frac{1}{1}}
  + x_{i}^{T} M_{i} u_{i} \right) +
  \frac{1}{2}x_N^T Q_N x_N + q_N^T x_N \\
& \mbox{s.t.} & x_0 = s_0 \\
& & x_{i+1} = A_i x_i + B_i u_i + c_i, \forall i \in \lbrace 0, \ldots, N-1 \rbrace.
\end{aligned}
\end{equation}
Note that the $R_i$ are required to be positive definite, and that the $Q_i - M_i R_i^{-1} M_i^T$
are required to be positive semi-definite, otherwise a minimum may not exist.

The conventional method for solving LQR problems sequentially relies on the fact that
the optimal cost-to-go functions are quadratic on the state at the corresponding stage.
Therefore, such functions $V_i(x_i) = 0.5 x_i^T P_i x_i + p_i^T x_i + z_i$ can be computed
in decreasing order of stage (i.e. $i$). Moreover, the optimal controls to be applied
at each stage can be shown to be affine functions of the state at the corresponding stage,
i.e. $u_i = K_i x_i + k_i$. The $K_i$ and $k_i$ can be computed as part of the same
backward pass as the $P_i$ and $p_i$.
Note that the constants $z_i$ do not need to be computed.
Initializing $P_N = Q_N$ and $p_N = q_N$, and introducing $G_i, H_i, h_i$ as convenient
auxiliary variables, the following recursion rules can be used:
\begin{equation}
\begin{aligned}
&G_i = R_i + B_i^T P_{i+1} B_i \\
&H_i = B_i^T + P_{i+1} A_i + M_i^T \\
&h_i = B_i^T (p_{i+1} + P_{i+1} c_i) + r_i \\
&K_i = -G_i^{-1} H_i \\
&k_i = -G_i^{-1} h_i \\
&P_i = Q_i + A_i^T P_{i+1} A_i + K_i^T H_i \\
&p_i = q_i + A_i^T (p_{i+1} + P_{i+1} c_i) + K_i^T h_i \\
\end{aligned}
\end{equation}
Once the $K_i, k_i$ have been computed, the $u_i$ can be computed in increasing order of $i$,
by alternating evaluations of $u_i = K_i x_i + k_i$ and $x_{i+1} = A_i x_i + B_i u_i + c_i$.
This is usually called the LQR forward pass.

\section{Associative Scans Overview}\label{associative-scans}

Associative scans are a common parallelization mechanism used in functional programming,
first introduced in~\cite{ref-AS}. They were used in~\cite{ref-AS-LQR} to derive a simple
method for solving (primal) LQR problems in $O(\log(m) + \log(N) \log(n))$ parallel time,
where $N, n, m$ are respectively the number of stages, states, and controls.

Given a set $\mathcal{X}$, a function $f: \mathcal{X} \times \mathcal{X} \rightarrow \mathcal{X}$
is said to be associative if $\forall a, b, c \in \mathcal{X}, f(f(a, b), c) = f(a, f(b, c))$.
The forward associative scan operation $S_f(x_1, \ldots, x_n; f)$ can be recursively defined by
$S_f(x_1; f) = x_1$ and $S_f(x_1, \ldots, x_{i + 1}; f) = (y_1, \ldots, y_i, f(y_i, x_{i + 1}))$,
where $y_1, \ldots, y_i = S_f(x_1, \ldots, x_i; f)$.
Similarly, the reverse associative scan operation $S_r(x_1, \ldots, x_n; f)$ can be recursively
defined by $S_r(x_1; f) = x_1$ and $S_r(x_1, \ldots, x_{i + 1}; f) = (f(x_1, y_2), y_2, \ldots, y_{i + 1})$,
where $y_2, \ldots, y_{i + 1} = S_r(x_2, \ldots, x_{i + 1}; f)$.
\cite{ref-AS} provides a method for performing associative scans of $N$ elements in parallel time $O(\log(N))$.

\section{Parallel Primal LQR Overview}\label{parallel-lqr}

A parallel algorithm for solving these problems was presented in~\cite{ref-AS-LQR},
although we require some minor changes due to a difference in problem formulations.
In order to describe this method, we need to introduce an important definition.
An interval value function $V_{i \rightarrow j}(x_i, x_j)$ maps states $x_i, x_j$ at stages $i, j$
to the minimum possible cost incurred in stages $i, \ldots, j-1$ among all trajectories that start
at state $x_i$ in stage $i$ and end at state $x_j$ in stage $j$ ($\infty$ if no such trajectory exists).
The key insight of~\cite{ref-AS-LQR} is that the optimal interval value functions can be represented as
\begin{equation}
\begin{aligned}
V_{i \rightarrow j} (x_i, x_j) = \max \limits_{\lambda} \left(
  \frac{1}{2} x_i^T P_{i \rightarrow j} x_i + p_{i \rightarrow j}^T x_i \right. \\
  \left. - \frac{1}{2} \lambda^T C_{i \rightarrow j} \lambda
  - \lambda^T \left( x_j - A_{i \rightarrow j} x_i - c_{i \rightarrow j} \right) \right).
\end{aligned}
\end{equation}
Initializing, for $i \in \lbrace 0, \ldots, N - 1 \rbrace$, $V_{i \rightarrow i+1}$ with
\begin{equation}
\begin{aligned}
  P_{i \rightarrow i+1} &= Q_i - M_i R_i^{-1} M_i^T, \\
  p_{i \rightarrow i+1} &= q_i - M_i R_i^{-1} r_i, \\
  A_{i \rightarrow i+1} &= A_i - B_i R_i^{-1} M_i^T, \\
  C_{i \rightarrow i+1} &= B_i R_i^{-1} B_i^T, \\
  c_{i \rightarrow i+1} &= c_i - B_i R_i^{-1} r_i, \\
\end{aligned}
\end{equation}
and initializing $V_{N \rightarrow N+1}$ with
\begin{equation}
\begin{aligned}
  P_{N \rightarrow N+1} &= Q_N, \\
  p_{N \rightarrow N+1} &= q_N, \\
  A_{N \rightarrow N+1} &= 0, \\
  C_{N \rightarrow N+1} &= 0, \\
  c_{N \rightarrow N+1} &= 0, \\
\end{aligned}
\end{equation}
the following combination rules can be applied to compute $V_{i \rightarrow k}$
from $V_{i \rightarrow j}$ and $V_{j \rightarrow k}$:
\begin{equation}
\begin{aligned}
  P_{i \rightarrow k} &= A_{i \rightarrow j}^T \left( I + P_{j \rightarrow k} C_{i \rightarrow j} \right)^{-1} P_{j \rightarrow k} A_{i \rightarrow j} + P_{i \rightarrow j}, \\
  p_{i \rightarrow k} &= A_{i \rightarrow j}^T \left( I + P_{j \rightarrow k} C_{i \rightarrow j} \right)^{-1} \left( p_{j \rightarrow k} + J_{j \rightarrow k} c_{i \rightarrow j} \right) \\
    & + p_{i \rightarrow j}, \\
  A_{i \rightarrow k} &= A_{j \rightarrow k} \left( I + C_{i \rightarrow j} P_{j \rightarrow k} \right)^{-1} A_{i \rightarrow j}, \\
  C_{i \rightarrow k} &= A_{j \rightarrow k} \left( I + C_{i \rightarrow j} P_{j \rightarrow k} \right)^{-1} C_{i \rightarrow j} A_{j \rightarrow k}^T + C_{j \rightarrow k}, \\
  c_{i \rightarrow k} &= A_{j \rightarrow k} \left( I + C_{i \rightarrow j} P_{j \rightarrow k} \right)^{-1} \left( c_{i \rightarrow j} - C_{i \rightarrow j} p_{j \rightarrow k} \right) \\
    & + c_{j \rightarrow k}. \\
\end{aligned}
\end{equation}
A reverse associative scan~\cite{ref-AS} can be used to compute the $P_{i \rightarrow N+1}, p_{i \rightarrow N+1}$
(i.e.~the $P_i, p_i$ from~\cref{sequential-lqr}),
which in turn can be used to compute the $K_i, k_i$ in $O(1)$ parallel time.

Finally, we wish to compute the $x_i, u_i$ from the $K_i, k_i$. Note that the sequential LQR forward pass
has $O(T)$ parallel time complexity. However, as done in~\cite{ref-AS-LQR}, we can reduce the computation of
the $x_i$ to a sequential composition of affine functions, which can also be parallelized with an associative scan~\cite{ref-AS}.
This will be described in the next section.

\begin{equation}
\begin{aligned}
  x_{i+1} &= A_i x_i + B_i u_i + c_i = A_i x_i + B_i \left( K_i x_i + k_i \right) + c_i \\
  &= \left( A_i + B_i K_i \right) x_i + \left( B_i k_i + c_i \right)
\end{aligned}
\end{equation}

Once these affine functions have been composed, they can be independently applied to $x_0$ to recover all the $x_i$.
The $u_i$ can then be computed in $O(1)$ parallel time by independently evaluating $u_i = K_i x_i + k_i$.

\section{Composing Affine Functions with Associative Scans}\label{composing-affine-functions}

In this section, we describe an $O(\log(N))$ parallel time algorithm for composing
$N$ affine functions $F_i(x) = M_i x_i + m_i$, as done in~\cite{ref-AS-LQR}.

Letting $\mathcal{X} = \mathbb{R}^{n} \times \mathbb{R}^{n \times n}$ and letting
$f: \mathcal{X} \rightarrow \mathcal{X}$ be defined by $f((a, B), (c, D)) = (Da + c, DB)$,
we claim that $f$ is associative. This is simple to verify:
\begin{equation}
\begin{aligned}
  &f(f((a, B), (c, D)), (e, F)) \\
= &f((Da + c, DB), (e, F)) \\
= &(F(Da + c) + e, F(DB)) \\
= &((FD)a + Fc + e, (FD)B) \\
= &f((a, B), (Fc + e, FD)) \\
= &f((a, B), f((c, D), (e, F))).
\end{aligned}
\end{equation}

Moreover, note that $f$ is the affine function composition operator, as
$D(Bx + a) + c = (DB)x + (Da + c)$.

\section{Dual LQR Backward Pass}\label{dual-lqr-backward-pass}

Next, we show how to compute the multipliers \(\lambda\),
i.e.~how to solve the dual part of the LQR problem.
Note that primal-dual LQR problems can be solved in a single Newton step,
independently of the starting values of $x, \lambda$. Therefore, starting
with $x, \lambda = 0$, the $\Delta x, \Delta \lambda$ from the Newton-KKT system
correspond to the optimal $x, \lambda$.
Having determined \(x\) using one of the primal LQR algorithms, determining $\lambda$
can be done by specializing
\(\nabla_{xx} \mathcal{L}(0) x + J(c)(0)^T \lambda =
-\nabla_x \mathcal{L}(0)\) (i.e.~the first block-row of the Newton-KKT system).
Note that we discard the rows where \(J(c)(x)^T\) contains any of the \(B_i\)
(as they are not useful for computing the \(\lambda_i\)),
and keep the rows where it contains any of the \(A_i\).

Specifically, for $i \in \lbrace 0, \ldots, N - 1 \rbrace$, it holds that
\begin{equation}
\begin{aligned}
&
Q_i x_i + M_i u_i
- \lambda_i + A_i^T \lambda_{i + 1} = -q_i
\Rightarrow \\
&
\lambda_i = Q_i x_i + M_i u_i
+ A_i^T \lambda_{i + 1} + q_i, \\
\end{aligned}
\end{equation}
and that
\begin{equation}
\begin{aligned}
&
Q_N x_N - \lambda_N = -q_N
\Rightarrow \\
&
\lambda_N = Q_N x_N + q_N.
\end{aligned}
\end{equation}
Note that these equations can be used to recursively compute
the \(\lambda_i\), in decreasing order of \(i\).

Moreover, note that these \(\lambda_i\) also satisfy the
block-rows of
\(\nabla_{xx} \mathcal{L}(0) x + J(c)(0)^T \lambda =
-\nabla_x \mathcal{L}(0)\) that we discarded, as a solution to
the full Newton-KKT system must exist due to the LICQ conditions being
satisfied, and as these \(\lambda_i\) are the \emph{only}
solution that satisfies the block-rows we did not discard.

\section{Parallel Dual LQR}\label{dual-lqr-associative-scan}

While the dual LQR backward pass described in~\cref{dual-lqr-backward-pass}
can be parallelized with a reverse associative scan (in a similar fashion to~\cref{composing-affine-functions}),
we can actually compute the multipliers $\lambda$ in $O(1)$ parallel time.

The multiplier $\lambda_0$ is also the Lagrange/KKT multiplier of
\begin{equation}
\min_{x_0} \frac{1}{2} x_0^T P_0 x_0 + p_0^T x_0 \quad \mbox{s.t.} \quad x_0 = s_0.
\end{equation}
By the first order optimality conditions, the gradient of the corresponding Lagrangian
\begin{equation}
\mathcal{L}(x_0, \lambda_0) = \frac{1}{2} x_0^T P_0 x_0 + p_0^T x_0 + \lambda^T (s_0 - x_0)
\end{equation}
must be $0$. In particular,
\begin{equation}
\nabla_{x_0} \mathcal{L}(x_0, \lambda_0) = 0 \Rightarrow \lambda_0 = P_0 x_0 + p_0.
\end{equation}

Similarly, each multiplier $\lambda_i$ is also the Lagrange/KKT multiplier of
\begin{equation}
\begin{aligned}
&\min_{u_i, x_{i+1}} &
\frac{1}{2} u_{i}^{T} R_{i} u_{i} + r_{i}^{T} u_{i} + x_{i}^{T} M_{i} u_{i} \\
  &&+ \frac{1}{2} x_{i+1}^T P_{i+1} x_{i+1} + p_{i+1}^T x_{i+1} \\
& \mbox{s.t.} & x_{i+1} = A_i x_i + B_i u_i + c_i.
\end{aligned}
\end{equation}
The corresponding Lagrangian is
\begin{equation}
\begin{aligned}
\mathcal{L}(u_i, x_{i+1}, \lambda_{i+1}) &= 
\frac{1}{2} u_{i}^{T} R_{i} u_{i} + r_{i}^{T} u_{i} + x_{i}^{T} M_{i} u_{i} \\
&+ \frac{1}{2} x_{i+1}^T P_{i+1} x_{i+1} + p_{i+1}^T x_{i+1} \\
&+ \lambda_{i+1}^T \left( A_i x_i + B_i u_i + c_i - x_{i+1} \right).
\end{aligned}
\end{equation}
Similarly, the first order optimality conditions imply
\begin{equation}
\nabla_{x_{i+1}} \mathcal{L}(u_i, x_{i+1}, \lambda_{i+1}) = 0 \Rightarrow
\lambda_{i+1} = P_{i+1} x_{i+1} + p_{i+1}.
\end{equation}

Finally, note that all the $\lambda_i$ can be evaluated independently,
so a backward pass is not required using these formulas.
This method of computing the multipliers also has the advantage of involving
fewer computations and being more numerically stable.

\section{Benchmarks}\label{benchmarks}

\begin{figure}[H]
\includegraphics[scale=0.25]{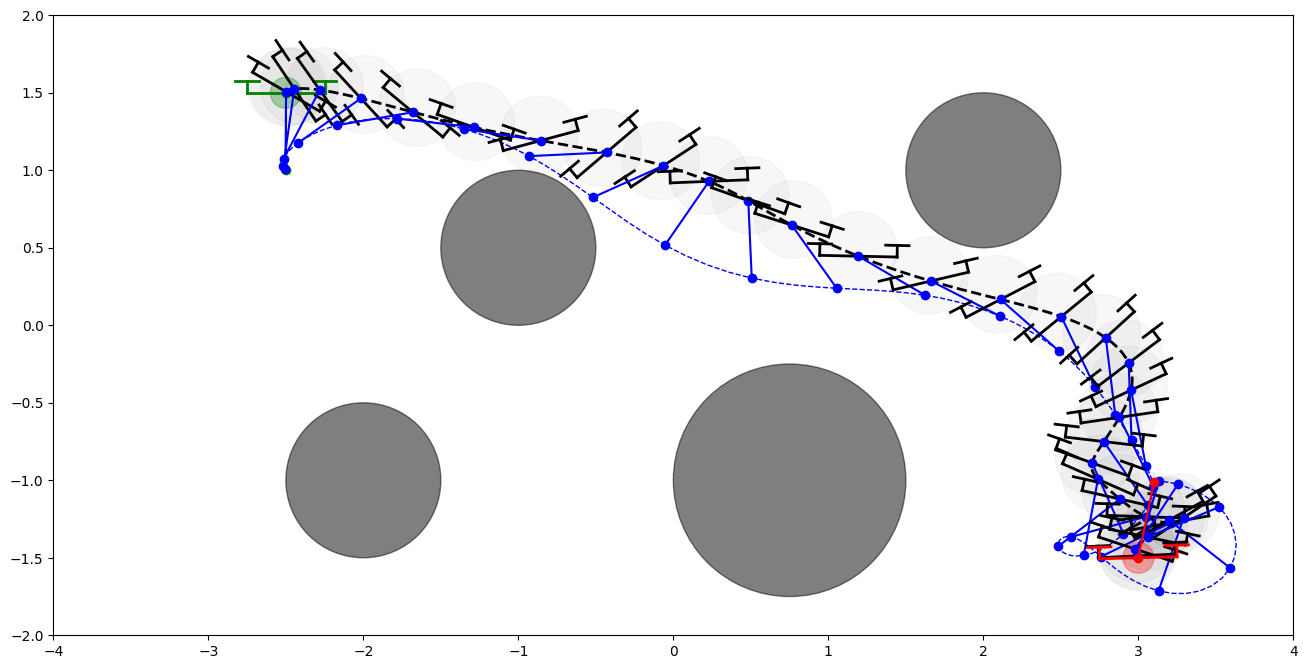}
  \caption{The solution of the quad-pendulum problem.}
  %A quadcopter with a pendulum avoiding obstacles to reach a goal state.
\end{figure}

Below, we include an example run of our algorithm on on a quad-pendulum
problem, which we got from \href{https://github.com/google/trajax}
{GitHub ($google/trajax$)} and modified
to treat constraints as high-penalty costs (using a gain of \(100\)). We
also increased the final state cost gains on the positions to \(1000\).
This seems to be the problem used in~\cite{ref-CLSQP}.
Note that we solve a single primal-dual LQR problem for each iteration,
and have to evaluate the user model once per line search step.
A JAX implementation of our algorithm takes 46ms to finish solving this problem
(on the CPU of an M2 MacBook Air).
Below, you can find a visual description of the problem and solution,
as well as solver logs.

\tablefirsthead{
  \toprule
  Iteration & Objective & $\|c\|^2$ & $m_\rho'$ & $\alpha$ \\
  \midrule}
\tablehead{%
\multicolumn{5}{c}%
% Uncomment for allowing the table to be split.
{{\bfseries  Continued from the previous column}} \\
% Comment above for allowing the table to be split.
\toprule
Iteration & Objective & $\|c\|^2$ & $m_\rho'$ & $\alpha$ \\ \midrule}
{
\tiny
% Uncomment for allowing the table to be split.
\tabletail{%
\midrule \multicolumn{5}{r}{{Continued on the next column}} \\ \midrule}
\tablelasttail{ \\\midrule }
% Comment above for allowing the table to be split.
\begin{supertabular}{lllll}
1  & 69.1912   & 9.173966407775879      & -196341.22  & 1.0 \\
2  & 38.890675 & 10.204395294189453     & -282.203    & 0.5 \\
3  & 21.515545 & 9.039937973022461      & -99.97078   & 0.25 \\
4  & 12.537336 & 3.5910232067108154     & -37.347027  & 0.5 \\
5  & 13.005388 & 1.256274938583374      & -7.008366   & 0.5 \\
6  & 11.193476 & 0.05883491784334183    & -8.02231    & 1.0 \\
7  & 10.798406 & 0.013367442414164543   & -0.79570323 & 1.0 \\
8  & 10.654287 & 0.0059180837124586105  & -0.2289545  & 1.0 \\
9  & 10.586153 & 0.001393172424286604   & -0.10068051 & 1.0 \\
10 & 10.561506 & 0.0004326050984673202  & -0.03127028 & 1.0 \\
11 & 10.545758 & 0.00018692569574341178 & -0.01689804 & 1.0 \\
12 & 10.535523 & 8.762026118347421e-05  & -0.01067685 & 1.0 \\
13 & 10.528571 & 4.414142313180491e-05  & -0.00720222 & 1.0 \\
14 & 10.523667 & 2.3347358364844695e-05 & -0.00506186 & 1.0 \\
15 & 10.520094 & 1.2880265785497613e-05 & -0.00367544 & 1.0 \\
16 & 10.517421 & 7.388789981632726e-06  & -0.0027444  & 1.0 \\
17 & 10.515375 & 4.351216375653166e-06  & -0.00209233 & 1.0 \\
18 & 10.513782 & 2.6609566248225747e-06 & -0.00163243 & 1.0 \\
19 & 10.512507 & 1.6749986571085174e-06 & -0.00129717 & 1.0 \\
20 & 10.511484 & 1.07657911030401e-06   & -0.00104529 & 1.0 \\
21 & 10.510648 & 7.084668709467223e-07  & -0.00085156 & 1.0 \\
22 & 10.509963 & 4.725115729797835e-07  & -0.00069964 & 1.0 \\
23 & 10.509394 & 3.217612913886114e-07  & -0.00057881 & 1.0 \\
24 & 10.50892  & 2.2059847992750292e-07 & -0.00048149 & 1.0 \\
25 & 10.508526 & 1.53070544683942e-07   & -0.00040247 & 1.0 \\
26 & 10.508192 & 1.0723919530164494e-07 & -0.00033782 & 1.0 \\
27 & 10.507912 & 7.579520655554006e-08  & -0.00028446 & 1.0 \\
28 & 10.507676 & 5.401408387228912e-08  & -0.0002403  & 1.0 \\
29 & 10.507476 & 3.848523988381203e-08  & -0.00020355 & 1.0 \\
30 & 10.507305 & 2.7843480410183474e-08 & -0.00017285 & 1.0 \\
31 & 10.50716  & 2.0009650825159042e-08 & -0.00014714 & 1.0 \\
32 & 10.507039 & 1.4696987626905411e-08 & -0.00012551 & 1.0 \\
33 & 10.506933 & 1.0591413968086272e-08 & -0.00010727 & 1.0 \\
\end{supertabular}%
}

\section{Conclusion}\label{conclusion}

This paper introduces a new algorithm for solving unconstrained discrete-time
optimal control problems, called Primal-Dual iLQR.

This framework can easily be extended to handle constrained problems.
Ideally, an interior point~\cite{ref-FATROP} method would be used to handle
general inequality constraints.
It would also be possible to handle both general equality and inequality constraints
with an augmented Lagrangian~\cite{ref-ALTRO} method, although local superlinear convergence
would be lost.
In both cases, the LQR shape of the posed subproblems can be preserved.
Preserving superlinar local convergence in the presence of general equality constraints
would require employing a linear equality-constrained LQR subproblem solver.
Using SQP~\cite{ref-ACADOS} to handle
general nonlinear constraints would yield Newton-KKT systems that are
fundamentally different from the ones we consider, and would therefore
require substantial changes to our method (including adding substantially
slower iterative subproblem solves).

Our method is substantially easier to warm start as compared to single-shooting
methods such as DDP, iLQR, or Stagewise Newton, as it treats both state and
control trajectories as free variables in the optimization problem (without
incurring extra computational costs); this allows them to be independently seeded,
even in dynamically infeasible ways.

The proposed algorithm is globally convergent and does not impede superlinear local
convergence, without requiring second order corrections to be applied.

Finally, solving the subproblems posed by our method can be parallelized to
run in time logarithmic in the number of stages, states, and controls.
Moreover, as our method never performs nonlinear dynamics rollouts,
it allows that the evaluation of the dynamics and derivatives
be computed in parallel, resulting in $O(1)$ parallel time per line search iteration.

\end{document}